\documentclass[12pt]{article}

\usepackage[colorlinks,allcolors=blue]{hyperref}

\usepackage[margin=1in]{geometry}
\usepackage{amsfonts, amsmath, amssymb, amsthm} 
\usepackage[none]{hyphenat}
\usepackage{fancyhdr} 
\usepackage[utf8]{inputenc}
\usepackage{mathtools}

\makeatletter
\def\legendre@dash#1#2{\hb@xt@#1{%
  \kern-#2\p@
  \cleaders\hbox{\kern.5\p@
    \vrule\@height.2\p@\@depth.2\p@\@width\p@
    \kern.5\p@}\hfil
  \kern-#2\p@
  }}
\def\@legendre#1#2#3#4#5{\mathopen{}\left(
  \sbox\z@{$\genfrac{}{}{0pt}{#1}{#3#4}{#3#5}$}%
  \dimen@=\wd\z@
  \kern-\p@\vcenter{\box0}\kern-\dimen@\vcenter{\legendre@dash\dimen@{#2}}\kern-\p@
  \right)\mathclose{}}
\newcommand\legendre[2]{\mathchoice
  {\@legendre{0}{1}{}{#1}{#2}}
  {\@legendre{1}{.5}{\vphantom{1}}{#1}{#2}}
  {\@legendre{2}{0}{\vphantom{1}}{#1}{#2}}
  {\@legendre{3}{0}{\vphantom{1}}{#1}{#2}}
}
\def\dlegendre{\@legendre{0}{1}{}}
\def\tlegendre{\@legendre{1}{0.5}{\vphantom{1}}}
\makeatother

\usepackage{booktabs, siunitx,caption}

\newtheorem{theorem}{Theorem}

\title{Effective Primality Test for $p2^n+1$, $p$ prime, $n>1$}
\author{Tejas R. Rao}
\date{}

\begin{document}
\maketitle

\begin{center}
\small \textbf{Abstract.} We develop a simple $O((\log n)^2)$ test as an extension of Proth's test for the primality for $p2^n+1$, $p>2^n$. This allows for the determination of large, non-Sierpinski primes $p$ and the smallest $n$ such that $p2^n+1$ is prime. If $p$ is a non-Sierpinski prime, then for all $n$ where $p2^n+1$ passes the initial test, $p2^n+1$ is prime with $3$ as a primitive root or is primover and divides the base $3$ Fermat Number, $GF(3,n-1)$. We determine the form the factors of any composite overpseudoprime that passes the initial test take by determining the form that factors of $GF(3,n-1)$ take. 
\end{center}

\setlength{\parindent}{0cm}

Proth's Theorem is well know throughout mathematics: where $P=k2^n+1$, $k$ odd, and $2^n>k$, 
\begin{center}
$a^{\frac{P-1}{2}}\equiv -1\mod P\Longleftrightarrow P$ prime, 
\end{center}
for all $a$ where $\legendre{a}{P}=-1$.
\vspace{5mm}\\
Sierpinski numbers of the second kind are integers $k$ such that $k2^n+1$ is not prime for all positive integers $n$ [\hyperlink{1}{1}]. For prime $k=p$, when $n=1$, one can use the Sophie-Germain Primality Test to determine whether $2p+1$ is prime. For $p<2^n$, one can use Proth's Theorem to determine whether $p2^n+1$ is prime. 
\vspace{5mm}\\
In this paper we extend a similar condition for $k=p$ prime and no constraints on the relative sizes of $2^n$ and $p$. For the rest of this paper, we will denote $R$ as $R=p2^n+1$, $n>1$. This bridges the gap between the Sophie-Germain Primality Test and Proth's Theorem, for $n>1$ where $p>2^n$. We utilize the definition of overpseudoprimes and primover numbers found in [\hyperlink{2}{2}]. Recognize that overpseudoprimes are a type of pseudoprime (always composite), whereas primover numbers are a type of probable prime (may be composite or prime). Also, $GF(3,z)=3^{2^{z}}+1$. 

\begin{theorem}
For all such $R$,
\begin{center}
$3^{\frac{R-1}{2}}\equiv -1\mod R\Longleftrightarrow R$ is prime or $R$ divides $GF(3,n-1)$ and is primover.
\end{center}
\end{theorem}
For necessity, assume $R$ is prime. We note that when $p=3$, Proth's theorem always applies and thus are theorem is also satisfied. We can also assume $3\nmid R$, because then the equivalence above would not hold. We thus have that $k\equiv\pm 1\mod 3$ and $2^n\equiv\pm 1\mod 3$. We cannot have $k\not\equiv 2^n\mod 3$, because then $3|R$. Therefore, 
\begin{align*}
R&\equiv 1+1\\
&\equiv -1\mod 3.
\end{align*}
This means that $\legendre{R}{3}=-1$. Additionally, we have $R\equiv 1\mod 4$. This means we can write $R=4m+1$, for $m\in\mathbb{N}$. Using quadratic reciprocity [\hyperlink{3}{3}], 
\begin{align*}
\legendre{3}{R}&=(-1)^{\frac{3-1}{2}\frac{R-1}{2}}\legendre{R}{3}\\
&=(-1)^{(1)\frac{4m}{2}}(-1)\\
&=(-1)^{2m+1}\\
&=-1.
\end{align*}
Since we assume $R$ is prime, we have $3^{\frac{R-1}{2}}\equiv -1\mod R$. If $R$ is composite but divides $3^{2^{n-1}}+1$, then the order of $3$ modulo $R$ is $2^n$ because its primitive index in $3^u+1$ is $2^{n-1}$ by [\hyperlink{4}{4}]. Since the order is $n$, it is clear that the conditions in the theorem are satisfied. 
\vspace{5mm}\\
For sufficiency, we assume $R$ is composite and does not divide $GF(3,n-1)$. Due to the conditions, we have 
 $3^{R-1}\equiv 1\mod R$. Therefore, the multiplicative order of $3$ mod $R$ divides $R-1=p2^n$ but does not divide $p2^{n-1}$ by the conditions specified. The multiplicative order is thus precisely $2^n$. Furthermore,  $3^{\frac{R-1}{2}}\equiv -1\mod R$ implies that all factors of $R$ divide $3^{\frac{R-1}{2}}+1$. This means that no factors divide $3^{\frac{R-1}{2}}-1$. But since the order of $R$ is $2^n$, all factors of $R$ must therefore share this order. Therefore all of the factors $f$ of $R$ have multiplicative order of $3$ modulo $f$ as $2^n$. Therefore we can write 
\begin{center}
$3^{2^n}-1=(3^{2^{n-1}}-1)(3^{2^{n-1}}+1)\equiv 0\mod R$.
\end{center}
and deduce that $3^{2^{n-1}}+1\equiv 0\mod R$. We arrive at a contradiction: a composite solution must divide $GF(3,n-1)$. Since all divisors of Fermat numbers are primover since they all have the same primitive index [\hyperlink{4}{4}], we prove the conditions. 
\vspace{5mm}\\
\textbf{Remark.} If Proth's theorem is satisfied and/or if $p>	\frac{3^{2^n}+1}{2}$ we know $R$ is prime iff it satisfies the above condition. 
\vspace{5mm}\\
\textbf{Remark.} We can additionally check that the number $R$ is not a Fermat factor of $GF(3,n-1)$ to prove that $R$ is prime after it passes the initial test. 
\vspace{5mm}\\
Next we will prove a minor theorem regarding the factors of overpseudoprimes that satisfy the test.
\begin{theorem}
For base $3$, all odd factors of $GF(3,n)$ are of the form $k2^{n+1}+1$, $3\nmid k$ and $k$ odd, or $3m2^{n+2}+1$, where $m$ is not necessarily odd.
\end{theorem}
When a factor $f$ of $GF(3,n)$ is a quadratic nonresidue base $3$, then it must have that 
\begin{center}
$3^{\frac{f-1}{2}}\equiv \legendre{3}{f}=-1\mod f$, 
\end{center} 
because all factors of generalized Fermat numbers are primover and therefore Euler-Jacobi pseudoprimes (REF). But since $3^{2^{n+2}}\equiv 1\mod f$, we know that $2^{n+2}\nmid\frac{f-1}{2}$. It follows that $f=2^{a}k+1$, where $k$ is odd and $a<2^{n+2}$. But since all factors of Fermat numbers are of the form $m2^{n+1}+1$ by (REF), where $m$ is not necessarily odd, we know that $f=2^{n+1}k+1$.\\
If $f$ is a quadratic residue base $3$, then we know that 
\begin{center}
$3^{\frac{f-1}{2}}\equiv \legendre{3}{f}=1\mod f$.
\end{center} 
By the same logic expressed above, $f=2^{n+2}m+1$, where $m$ is not necessarily odd.\\
But since we can write $f=4m+1$, for $m\in\mathbb{N}$. Using quadratic reciprocity (REF), 
\begin{align*}
\legendre{3}{f}&=(-1)^{\frac{3-1}{2}\frac{f-1}{2}}\legendre{f}{3},
\end{align*}
and $\legendre{f}{3}\equiv f^{\frac{3-1}{2}}=f\mod 3$, the result follows. 
\vspace{5mm}\\
\textbf{Corollary.}
All factors $f$ of a composite $R=p2^n+1$, $p>3$, $n>2$, and $3^{\frac{R-1}{2}}\equiv -1\mod R$, are of the form
\begin{center}
$f=k2^n+1$, $k\in\mathbb{N}$,
\end{center}
and furthermore are primover with $O_f(3)=2^n$. 

\newcommand{\noopsort}[1]{} \newcommand{\printfirst}[2]{#1}
  \newcommand{\singleletter}[1]{#1} \newcommand{\switchargs}[2]{#2#1}


\begin{thebibliography}{10}

\bibitem{article-full}
Keller, W.\\
\textbf{Factors of Fermat Numbers and Large Primes of the Form $k2^n+1$}.\\
\emph{\href{https://cs.uwaterloo.ca/journals/JIS/VOL15/Castillo/castillo2.pdf}{Mathematics of Computation},} 41 (1983), 661-673 \hypertarget{1}{}

\bibitem{article-full}
Castillo, J. H., García-Pulgarín, G., Shevelev, V., Velásquez-Soto, J. M.\\
\textbf{Overpseudoprimes, and Mersenne and Fermat Numbers as Primover Numbers}.\\
\emph{\href{https://cs.uwaterloo.ca/journals/JIS/VOL15/Castillo/castillo2.pdf}{Journal of Integer Sequences},} 15 (2012), Issue 7 \hypertarget{2}{}

\bibitem{article-full}
Beiler, A. H.\\
\textbf{Recreations in Theory of Numbers.}\\
\emph{Dover Publications, LCCN \href{https://catalog.loc.gov/vwebv/holdingsInfo?searchId=3101&recCount=25&recPointer=0&bibId=12849868}{64013458},} (1964) \hypertarget{3}{} 

\bibitem{article-full}
Rao, T. R.\\
\textbf{Primitive Indexes, Zsigmondy Numbers, and Primoverization}.\\
\emph{\href{https://arxiv.org/abs/1810.11456}{Cornell University Library},} (2018) \hypertarget{4}{}

\end{thebibliography}
\end{document}